\documentclass[11pt]{article}

\usepackage{a4}
\usepackage{latexsym}
\usepackage{amsthm}
\usepackage{amsmath}
\usepackage{amsfonts}
\usepackage{german}

\usepackage{graphics}

\theoremstyle{plain}
\newtheorem{thm}{Theorem}[section]
\newtheorem*{ToePC}{T\"oplitz Pencil Conjecture (ToePC)\footnote{we cannot abbreviate TPC, for this denotes the twin primes conjecture.} (\cite{ScSh1})}

\newtheorem{observation}[thm]{Observation}
\newtheorem{property}[thm]{Property}
\newtheorem{criterion}[thm]{Singularity criterion }

\theoremstyle{definition}

\newcommand{\tr}{{}^t}

\newcommand{\C}{{\mathbb C}}

\newcommand{\bmx}{\left[\begin{matrix}}
\newcommand{\emx}{\end{matrix}\right]}
\newcommand{\Qn}{Q^{(n)} }
\newcommand{\iQn}{Q^{(n)^{-1}} }

\title{Algebraic condition for the singularity of certain T\"oplitz pencils}
\author{Wiland Schmale\thanks{ University of Oldenburg, Germany,  email: wiland.schmale@uni-oldenburg.de}}
\date{ June 21, 2017}

\begin{document}
\parindent3mm	
\maketitle
\begin{abstract}
An algebraic condition for the singularity of certain T\"oplitz matrix pencils is derived which involves only the principal minors of the constant parts of the pencils. This leads to an algebraic conjecture which is equivalent to the so-called T\"oplitz pencil conjecture.
\end{abstract}

keywords: T\"oplitz pencil, conjecture, singularity\ \ \ 
MSC2010: 15B05, 93C05

\section{Introduction}
Let $T(x)=M_0+xM_1$, where the $n\times n$-matrices $M_0, M_1$ are given as

\[
M_0={\small\bmx c_2&c_1&0&\dots&0\\
          c_3&c_2&c_1&\dots&0 \\
          \vdots&\vdots && \ddots&\vdots \\
          c_{n}&c_{n-1}& \dots& \dots&c_1\\
          c_{n+1}&c_n& \dots& \dots&c_2\\
          \emx},
\ M_1={\small\bmx 0&0&1&0&\dots&0\\
          0&0&0&1&\dots&0 \\
          &\vdots&\vdots && \ddots&\vdots \\
          0&0& \dots & \dots &\dots&1\\
          0&0& \dots& \dots&\dots&0\\
          0&0& \dots& \dots&\dots&0\\
          \emx},
           n\ge 2.
\]          
The complex numbers $c_1,c_2,\dots,c_{n+1}$ are all supposed to be 
{\bf non-zero}. $x$ is an indeterminate over $\C$. $M_1$ is the $2\times 2$ zero-matrix if $n=2$.

\begin{ToePC}
If $\det T(x)=0$ as a polynomial in $\C[x]$, i.e. all coefficients vanish, then the first two columns of $M_0$ (or $T(x)$) are linearly dependent. I.e. for some non zero complex number $\lambda$ one has $c_{k+1}=\lambda c_k$  or equivalently $c_{k+1}=\lambda^kc_1$ for $k=1,\dots,n$. 
\end{ToePC}

The converse of the conjecture is trivial. Some authors prefer an equivalent Hankel version e.g. \cite{GoKo}.

The conjecture originates from an older (1981) and still open conjecture from linear systems theory (\cite{BSSV}, p. 124) on feedback cyclization over $\C[x]$. More background information and references can be found e.g. in \cite{ScSh1}. 

To date despite several efforts ToePC is still open (\cite{ScSh2,BoLi,HW,open}). The best result till now is achieved in \cite{GoKo}, where its truth could be established for $n\le 8$. 
A proof of the truth of ToePC would at the same time be a further progress in the proof of the older conjecture from linear system theory.\\

The various approaches exploit the singularity of $T(x)$ in different ways. 
In these notes, based on an old observation for singular matrix pencils and related block matrices, a new necessary and sufficient algebraic condition for the singularity of the specific pencils $T(x)$ occurring in ToePC will be derived. A closer look at this condition then will lead to a condition relying exclusively on the principal minors of $M_0$. At the same time this leads to a new algebraic conjecture equivalent to ToePC.\\

The polynomial $\det T(x)$ is of degree $n-2$ and its $k$-th coefficient is homogeneous in $c_1,\dots,c_{n+1}$ of degree $n-k$.
Therefore ToePC is true if and only if it is true, say with $c_1=1$. This will only be exploited in Section 3 to simplify expressions. In \cite{GoKo} by tricky arguments it is shown that without loss one can even assume that $c_1=1$ and $c_2=1$. 

\section{From block matrices to matrix powers}
If $\det T(x)=0$, then $T(x)$ is a singular matrix pencil. Singularity implies the existence of a non zero vector polynomial $f(x)$ satisfying $T(x)f(x)=0$. {\bf Let \pmb{$d$} denote the minimal degree of such vector solutions.} An old observation from the theory of singular pencils (see e.g. \cite{G}, p. 29) then states:\\

\begin{observation} \label{oO}
  Let {$\det T(x)=0$} and {$d\ge 0$}.
  
\noindent   The equation 
{$T(x)f(x)=0$} admits only solutions of degree {$\ge d$}, if and only if  the matrix 
\[
{
\mathcal{C}=\small
\underbrace{
\bmx
     M_0&0&&0&0&0\\
     M_1&M_0&&0&0&0\\
     0&M_1&\rotatebox{-20}{$\ddots$}&0&0&0\\
     \vdots&\vdots&\rotatebox{-20}{$\ddots$}&&\vdots&\vdots\\    
     0&0&&M_0&0&0\\
     0&0&\cdots&M_1&M_0&0\\
     0&0&\cdots&0&M_1&M_0\\
     0&0&\cdots&0&0&M_1
\emx
}_{\displaystyle d+1 \text{\ block columns}}
}
\]
has the property:

the first {$n\cdot d$} columns are linearly independent (l.i.) and at the same time all columns ($n\cdot (d+1)$) are linearly dependent (l.d.).
 The result is trivial for $d=0$.
\end{observation}
\noindent The if-part of the observation is not explicitely stated in \cite{G} but is immediate.
For later use $\mathcal{C}$ is displayed for some $d\ge 5$. For $d=0$, $d=1$ etc.  one has 
$\mathcal{C}=\small\bmx M_0\\M_1\emx$, $\mathcal{C}=\small\bmx M_0&0\\M_1&M_0\\0&M_1\emx$, etc.

Starting from this observation we will derive now an algebraic necessary and sufficient condition for the singularity of $T(x)$. For this we subdivide the matrices $M_0,M_1$ as follows:
\[
	M_0={\small\bmx
			v&Q\\
			c_{n+1}&w
		\emx},
	M_1={\small\bmx
			0&B\\
			0&0
		\emx},
\]	
where
\[
	Q={\small\bmx
		c_1&\cdots&0\\
		\vdots&\ddots&\vdots\\
		c_{n-1}&\cdots&c_1
	  \emx},
	  v={\small\bmx
	  		c_2\\
	  		\vdots\\
	  		c_n
	    \emx},
	  w={\small\bmx
	  		c_n,&\dots&,c_2
	    \emx},
	  B={\small\bmx
	  		0&1&\cdots&0\\
	  		\vdots&&\ddots&\vdots\\
	  		0&0&&1\\
	  		0&0&\cdots&0
	  	\emx}.
\]
Note that  $w=\tr v P_{n-1}$, where 
$P_{n-1}=\small\bmx 0&\cdots&1\\\vdots&\ddots&\vdots\\1&\cdots&0\emx$
is the $(n-1)\times (n-1)$ - anti-diagonal unit-matrix and $\tr v$ the transpose of $v$.
Since $c_1\ne 0$ the matrix $Q$ is invertible.
According to our subdivision of the matrices $M_0,M_1$ the overall matrix $\mathcal C$ looks as follows:
\[
\mathcal C=\small
\underbrace{	
		\left[\begin{array}{cc|cc|cc|cc|cc|cc}
			v&Q\\
			c_{n+1}&w\\
			0&B&v&Q\\			
			0&0&c_{n+1}&w\\
			&&0&B\\
			&&0&0\\
			&&&&\ddots\\
			&&&&&\ddots\\			
			&&&&&&v&Q\\
			&&&&&&c_{n+1}&w\\									
			&&&&&&0&B&v&Q&\\
			&&&&&&0&0&c_{n+1}&w\\
			&&&&&&&&0&B&v&Q\\
			&&&&&&&&0&0&c_{n+1}&w\\
			&&&&&&&&&&0&B\\
			&&&&&&&&&&0&0\\
		\end{array}\right]
}_{\displaystyle d+1 \text{\ vertical blocks}}
\]

The following condition for the singularity of $T(x)$ can now be derived:\\

\begin{criterion}

$T(x)$ is singular if and only if 
\[
	\pmb{wQ^{-1}(BQ^{-1})^kv=0\ \text{ for all } k\ge 1 \text{ and } wQ^{-1}v=c_{n+1}\ (\ne 0)}\hspace{3cm}		\hfill {\pmb{(\mathcal{S})}}
\]
\end{criterion}

\begin{proof} \underline{First we derive the condition $\mathcal{(S)}$.} $Q$ can be used  for \glqq cleaning up\grqq\ in the matrix $\mathcal C$. In order to maintain the l.d./l.i.-properties in observation \ref{oO}, we will {\it not perform operations involving columns of the last block and the preceding ones at the same time.}

We first transform the $v$-entries to $0$ by column operations  within each block column and then the $w$- and $B$-entries to zero by row operations.
The resulting matrix will be $\mathcal{C'}$ below. In order to render a reasonable display  we introduce the abbreviation 
\begin{equation} \label{A}
	A:=(-BQ^{-1})
\end{equation}
The last three block columns indicate what $\mathcal{C'}$ looks like in case  $d=2$.

\[	
	\mathcal{C'}=\small\underbrace{
		\left[\hspace{-3mm}\begin{array}{cc|cc|cc|cc|cc|ccc}
				0&Q   \\
				\star&0&   \\
	 			Av&0&0&Q&  \\
				-wQ^{-1}Av&0&\star&0&\phantom{\ddots}\ddots   \\
				A^2v&0&Av&0&\phantom{\ddots}\ddots\\
				-wQ^{-1}A^2v&0&-wQ^{-1}Av&0 &\phantom{\ddots}\ddots   \\							
				\vdots&\vdots&\vdots&\vdots   \\				
				A^{d-2}v&0&A^{d-3}v&0&\phantom{\cdots}\cdots&& 0&Q\\
				-wQ^{-1}A^{d-2}v&0&-wQ^{-1}A^{d-3}v&0&\phantom{\cdots}\cdots&&\star&0\\
				A^{d-1}v&0&A^{d-2}v&0&\phantom{\cdots}\cdots&&Av&0&  0&Q\\
				-wQ^{-1}A^{d-1}v&0&-wQ^{-1}A^{d-2}v&0&\phantom{\cdots}\cdots&& -wQ^{-1}Av&0&  \star&0&&&\\
	 			A^{d}v&0&A^{d-1}v&0&\phantom{\cdots}\cdots&&A^2v &0&  Av & 0&0&Q\\
				-wQ^{-1}A^dv&0&-wQ^{-1}A^{d-1}v&0&\phantom{\cdots}\cdots&& -wQ^{-1}A^2v &0& -wQ^{-1}Av&0&\star&0\\
				A^{d+1}v&0&A^{d}v&0&\phantom{\cdots}\cdots&& A^3v&0& A^2v & 0&Av&0\\
				0&0&0&0&\phantom{\cdots}\cdots&& 0&0&  0&0&0&0\\
	\end{array}\hspace{-3mm}\right]	}_{\displaystyle d+1 \text{\ vertical blocks}}
\]
where $\star=c_{n+1}-wQ^{-1}v$.

Since $\det(T(x))=0$ in $\C[x]$ we must have $\det(M_0)=0$ and therefore the $\star$-entries in the matrix $\mathcal{C'}$ are zero. The latter means that $c_{n+1}=wQ^{-1}v$.  Note that $A^kv=(-BQ^{-1})^kv$ is a $(n-1)$-vector and  
$-wQ^{-1}A^kv=-wQ^{-1}(-BQ^{-1})^kv$ is a complex number.

We can further clean up the rows left of the $Q$-entries within the first $d$ block columns and then cancel the first $d$ columns with $Q$-entries. We also can cancel zero-rows. The remaining matrix with four items highlighted for later use is $\mathcal C''$:\newline
\[
	\mathcal{C''}=\small
\underbrace{\left[\begin{array}{ccccc|cc}	
			-wQ^{-1}Av\\			
			-wQ^{-1}A^2v&-wQ^{-1}Av\\
			\vdots&\vdots&\ddots\\
			-wQ^{-1}A^{d-2}v&-wQ^{-1}A^{d-3}v   \\
			-wQ^{-1}A^{d-1}v&-wQ^{-1}A^{d-2}v&\cdots&-wQ^{-1}Av   \\
			A^dv&A^{d-1}v&\cdots& \pmb{A^2v}   &\pmb{Av}&0&Q\\
			-wQ^{-1}A^dv&-wQ^{-1}A^{d-1}v&\cdots&-wQ^{-1}A^2v   &-wQ^{-1}Av&0&0\\
			A^{d+1}v&A^dv&\cdots&\pmb{A^3v}    &\pmb{A^2v}&Av&0\\			
	\end{array}\right]
}_{\displaystyle d \text{ columns\ } \text{ plus } 1 \text{ column } \text{ plus } (n-1) \text{ columns from } Q }
\]

 $\mathcal C''$ inherits from $\mathcal C$ the property:
 \begin{center}
	the first $d$ columns are l.i. and all columns together $(d+1+(n-1))$ are l.d.
\end{center}

As far as we only consider the linear dependence of all columns of $\mathcal C''$ we can also omit the third last row and last block column which gives $\mathcal{C'''}$:

{
\[
\mathcal{C'''}=\small
\underbrace{
	\left[\begin{array}{cccccccc}	
			-wQ^{-1}Av&\\			
			-wQ^{-1}A^2v&-wQ^{-1}Av\\
			\vdots&&\hspace{-3mm}\ddots\\
			-wQ^{-1}A^{d-2}v&-wQ^{-1}A^{d-3}v   \\
			-wQ^{-1}A^{d-1}v&-wQ^{-1}A^{d-2}v&\cdots&-wQ^{-1}Av\\			
			-wQ^{-1}A^dv&-wQ^{-1}A^{d-1}v&\cdots&-wQ^{-1}A^2v&-wQ^{-1}Av\\
			A^{d+1}v&A^{d}v&\cdots&A^3v&A^2v&Av\\			
	\end{array}\right]
}_{\displaystyle d+1 \text{ columns}}
\]
}

Now, if $Av=0$ then $\mathcal{C'''}$ is the $0$-matrix and $(\mathcal{S})$ follows trivially.\\

Let therefore {$Av\ne 0$}.
By linear dependence we conclude from $\mathcal C'''$: 
\[
	{wQ^{-1}Av=0.}
\]
If $wQ^{-1}A^2v\ne 0$, then (look at $\mathcal C'''$) the vectors 
$A^2v$ and $Av$ must be l.d., say 
$\rho\, A^2v+\sigma\, Av=0$. Therein $\rho$ must be non-zero since $Av\ne 0$. I.e. 
$A^2v=-\sigma/\rho\, Av$. Multiplication of the dependence relation from the left by $wQ^{-1}$ shows that $wQ^{-1}A^2v$ cannot be non-zero. So also
\[
	{wQ^{-1}A^2v=0}.
\]

Now assume $wQ^{-1}A^3v\ne 0$. Then once more by linear dependence  we must have a dependence relation between 
$A^3v,\ A^2v$ and $ Av$,
while in $\mathcal C''$ the submatrix
\[\small
	\bmx
		A^2v&Av\\
		0&0\\
		A^3v&A^2v
	\emx
\]
must have full rank. If $A^2v$ and $BQ^{-1}v$ are l.d. then also $A^3v$ and $A^2v$ with the same dependence relations. Therefore the rank can no longer be full. So the vectors are actually l.i. and in a non-trivial dependence relation for $A^3v,\ A^2v$ and $ Av$ the coefficient of $A^3v$ must be non-zero. Multiplying from the left by $wQ^{-1}$ shows finally contrary to our assumption that also {$wQ^{-1}A^3v=0$}.\\

The latter argument can be repeated inductively. Assume \mbox{$wQ^{-1}A^kv=0$} for $1\le k\le r$ for some $3\le r\le d-1$. In this case the matrix $\mathcal C''$ looks as follows:

\[
	\mathcal{C''}=\small
\underbrace{\left[\begin{array}{ccccccccc}	
			0&0\\			
			0&0\\
			\vdots&\vdots\\
			0&0\\
			\underline{-wQ^{-1}A^{{r+1}}v}&0\\[1.5mm]			
			-wQ^{-1}A^{{r+2}}v&\underline{-wQ^{-1}A^{{r+1}}v}\\
			\vdots&&\rotatebox{-25}{$\ddots$}&\\
			-wQ^{-1}A^{d-1}v&-wQ^{-1}A^{{d-2}}v&&0&0&&0&0&0\\
			A^dv&A^{d-1}v& &{A^{r+1}v}& \pmb{A^rv}&\cdots&\pmb{Av}&0&Q\\
			-wQ^{-1}A^dv&-wQ^{-1}A^{d-1}v&&\underline{-wQ^{-1}A^{r+1}v}&0&\cdots&0&0&0\\[1.5mm]
			A^{d+1}v&A^dv&\cdots&{A^{r+2}v}&\pmb{A^{r+1}v}&\cdots&\pmb{A^2v}&Av&0\\			
	\end{array}\right]
}_{\displaystyle (d+1) \text{ columns\ } + (n-1) \text{ columns from } Q}
\]

Now if $-wQ^{-1}A^{r+1}v\ne 0$ then the vectors 
\[{A^{r+1}v},\dots {A^2v},Av\] must be l.d. At the same time the submatrix 
\[\small
	\bmx
		{A^rv}&\cdots&{Av}\\
		0&\cdots&0\\
		{A^{r+1}v}&\cdots&{A^2v}
	\emx	
\]		
must have full rank which is only possible if the vectors in the first row are l.i. Therefore the vector $A^{r+1}v$ is in the span of the latter vectors. Multiplying from the left by $wQ^{-1}$ shows then that 
{$-wQ^{-1}A^{r+1}v=0$}. We conclude

\begin{equation}\label{eq1}
	wQ^{-1}A^kv=0 \text{ for } k=1,\dots,d
\end{equation}

and for the only remaining nonzero third last and last row we must have:
\begin{equation}\label{eq2}
	\text{The vectors } A^kv \text{ for } k=1,\dots,d \text{ are l.i. }
\end{equation}
\begin{equation}	\label{eq3}
	\text{and the vectors } A^kv \text{ for } k=1,\dots,d+1 \text{ are l.d. }
\end{equation}
From  (\ref{eq3}) and (\ref{eq2}) we conclude 
\begin{equation}\label{eq4}
A^{d+1}v=\sum_{k=1}^d\lambda_k A^kv 
\end{equation}

Multiplying by $wQ^{-1}$ gives 
$wQ^{-1}A^{d+1}v=0$. Multiplying successively (\ref{eq4}) by 
$A^s,\ s\ge 1 $ and then by $wQ^{-1}$, together with (\ref{eq1}) and remembering (\ref{A}), gives condition ($\mathcal S$).\\

\noindent\underline{ Now we assume  ${\mathcal{(S)}}$ and have to derive the singularity of ${T(x)}$:}\\[2mm]
Suppose $BQ^{-1}v=0$. 
The matrix $\small \bmx 0&Q\\0&0&\\-BQ^{-1}v&0\\0&0\emx$ is then singular.
By operations as executed above we can obtain the matrix
$\small\bmx
	v&Q\\
	wQ^{-1}v&w\\
	0&B\\
	0&0\\
\emx$, 
which still must be singular. Since $c_{n+1}=wQ^{-1}v$ the last matrix is just
$\small\bmx M_0\\M_1 \emx$. It is singular only if the first two columns of $M_0$ are dependent (over $\C$). In this case also $T(x)$ is singular over $\C[x]$.\\

Now we must assume $BQ^{-1}v\ne 0$. Let $d$ be minimal such that the $d+1$ vectors $A^{d+1}v$, $\dots$, $Av$ are l.d. Remember (\ref{A}). We can form the matrix $\mathcal {C'}$ as above, most entries being 0 by condition
($\mathcal S$) but with the same l.i/l.d-property. We then transform exactly as above but in the opposite direction to obtain $\mathcal C$ with l.i/l.d. properties as in the theorem. Therefore observation \ref{oO} tells us that $T(x)$ is singular admitting a non zero solution of minimal degree d for the equation $T(x)f(x)=0$.
\end{proof}

\section{Transition to principal minors of ${\pmb{M_0}}$}

It is worthwhile to have a closer look at the conditions $(\mathcal S)$. It turns out that $(\mathcal S)$ translates into an equivalent condition only for the principal minors of $M_0$. To avoid many powers of $c_1$ in the following we assume $c_1=1$. This is no restriction for ToepC as was mentioned at the beginning. The principal minors of $M_0$ will be denoted as  
$m_{r}\,, 1\le r\le n$, for example $m_1=c_2, m_{2}=c_2^2-c_3$,  etc. We set $m_0=1$.\\

One reason for the appearance of principal minors in our situation is that one has the relations
(assume just for the moment that  $c_2,\dots,c_{n+1}$ are independent variables over $\C$.)
\[
	m_{r}=(-1)^{r+1}c_{r+1}+p_r(c_2,\dots,c_r)
\]
for some polynomial expressions $p_r\in \C[c_2,\dots,c_{n+1}]$
and therefore 
\begin{equation}\label{var}
\C[c_2,\dots,c_{n+1}]=\C[m_1,\dots, m_{n}]
\end{equation}
The latter means that in this case also the minors can be considered as  independent variables.

Furthermore, if for some complex numbers $c_2,\dots,c_{n+1}$ we have $m_{2}=0,\dots,m_{n}=0$, then $c_k=c_2^{k-1}$ for $2\le k\le n+1$ and vice versa, both conditions beeing equivalent to the linear dependence of the first two columns of $M_0$ or $T(x)$.\\

The following properties are basic for our transition to minors. They seem not to be present explicitely in published material. Partly they appear in  unpublished notes of W.F. Trench (\cite{T}) in a different context but providing a hint for a  short proof of Property \ref{pr1} below. In order to indicate the size of $Q$, if desirable in induction steps, we write $\Qn$ instead of the $(n-1)\times (n-1)$ - matrix $Q$ and also ${M_0}^{(n)}$ instead of the $n\times n$ - matrix $M_0$. Also $I_n$ will denote the $n\times n$ - identity matrix.

\begin{property}\label{pr1}
$\iQn$ is also a lower triangular T\"oplitz matrix and 
\begin{equation}\label{epr1}
	\iQn={\small\bmx 1   &0&\cdots & &&0\\
				 -m_1&1&\cdots&&& 0 \\
				 m_{2}&-m_1&\cdots &&&0\\
				 -m_{3}&m_{2}&\cdots&&&0\\
				 \vdots&\vdots&\ddots&&&\vdots\\
				 \text{\small{$(-1)^{n-2} m_{n-2}$}}
				      &\text{\small{$(-1)^{n-3}m_{n-3}$}}&\cdots&-m_1&&1
			\emx}
			=\big( (-1)^{i+j}m_{i-j} \big)_{1\le i,j\le n-1},
\end{equation}
where we set $m_{i-j}=0$ if $i<j$. 
\end{property}

\begin{proof} 
$\iQn=\bmx 1\emx$ in case $n=2$. Suppose $n\ge 1$. One observes that
\[
	{Q^{(n+1)}}^{-1}=
	\bmx		
		\iQn&0\\
		\star\cdots&1
	\emx
	=
	\bmx
		1&\rotatebox{90}{$0$}\\
		\stackrel{\vdots}{\star}&\stackrel{\displaystyle\iQn}{\phantom{...}}
	\emx							
\] 
Because of the two overlapping blocks $\iQn$ by induction hypothesis ${Q^{(n+1)}}^{-1}$ is automatically T\"oplitz and we only have do determine the 
$[n,1]$-entry indicated by  $\star$. The classical adjoint matrix of $Q^{(n+1)}$ (see e.g. \cite{HoJo}, p. 22) tells us that
\[
	\big({Q^{(n+1)}}^{-1}\big)\pmb{[n,1]}=(-1)^{n+1}\cdot \det\underset{\pmb{\{1,n\}}}{Q}^{(n+1)}
\]

where $\underset{\{1,n\}}{Q}^{(n+1)}$ is the matrix $Q^{(n+1)}$ with row $1$ and column $n$ deleted. But this matrix is nothing else but $M_0^{(n-1)}$. Therefore the $[n,1]$ - entry of $({Q^{(n+1)}})^{-1}$ is
\[
	(-1)^{n+1}\cdot \det\underset{{\{1,n\}}}{Q}^{(n+1)}
	=(-1)^{n+1}\cdot \det M_0^{(n-1)}=(-1)^{n+1}\cdot m_{n-1}=(-1)^{n-1}\cdot m_{n-1}
\]
as claimed  in (\ref{epr1}), where $n$ is to be replaced by $(n+1)$.
\end{proof}

\begin{property}
\[
	\iQn v={\small
	\bmx
		m_1\\-m_{2}\\m_{3}\\\vdots\\
		\text{\small{$(-1)^{n-2}m_{n-1}$}}		
	\emx} \text{\ \ and as a consequence \ \ }		
	B\iQn v=\bmx\small
		-m_{2}\\m_{3}\\\vdots\\
		\text{\small{$(-1)^{n-2}m_{n-1}$}}\\0		
	\emx	
\]
\end{property}
\begin{proof}
In the equation ${Q^{(n+1)}}^{-1} \cdot Q^{(n+1)}=I_n$ we consider the first column of $I_n$ and recall property \ref{pr1}. This gives
\[{\small
	\bmx1\\0\\ \vdots\\0\emx}
	=\small{Q^{(n+1)}}^{-1}\cdot \bmx 1\\c_2\\\vdots\\c_n\emx=	
	\left[ \ \begin{array}{c|ccc}  
		1&&0\\
		\hline		
		-m_1&&\\
		\vdots&&\iQn\\
		(-1)^{n-1}m_{n-1}&
	\end{array}\  \right]\cdot \bmx 1\\c_2\\\vdots\\c_n\emx.
\]
From this equation the property can be read out. Note that one can choose any $c_{n+2}$ to obtain a matrix $M_0^{(n+1)}$ as a frame for $Q^{(n+1)}$ without affecting the proof.

\end{proof}

Let in addition 
\[
	X={\small
			\bmx 
				 -m_1&1&\cdots& 0 \\
				 m_{2}&-m_1&\cdots &0\\
				 -m_{3}&m_{2}&\cdots&0\\
				 \vdots&\vdots&\ddots&\vdots\\
				 \cdot&\cdot&&1\\
				 \text{\small{$(-1)^{n-2} m_{n-2}$}}
				      &\text{\small{$(-1)^{n-3}m_{n-3}$}}&\cdots&-m_1&
			\emx}
		\ \text{and}\
	y=\small\bmx
		-m_{2}\\m_{3}\\\vdots\\
		\text{\small{$(-1)^{n-2}m_{n-1}$}}		
	\emx\ .
\]
$X$ is obtained by cancelling the first row and the last column in $\iQn$.\\

\begin{property}
$\det X=\det X^{(n)}=(-1)^{n}\,c_{n-1}$, i.e. non-zero by assumption.
\end{property}

\begin{proof}
Once more the classical adjoint matrix (notation as in the proof of Property \ref{pr1}) tells us  that
\begin{eqnarray*}
	c_{n-1}	= \Qn[n-1,1]
	= {\big(\iQn\big)^{-1}}{[n-1,1]}\\
	= (-1)^{n}\cdot \det\underset{\{1,n-1\}}{\big(Q^{(n)}\big)}^{-1}
	= (-1)^{n}\cdot\det X^{(n)}
\end{eqnarray*}

\end{proof}

With the help of the matrix $X$ and the vector $y$ the conditions $(\mathcal S)$ translate into the following conditions ${(\mathcal{SM})}$ that involve {\bf only} the principal minors of $M_0$.

\begin{criterion}
$T(x)$ is singular if and only if ($X^0=I_{n-2}$)
\[
\pmb{(\tr yP_{n-2})\ X^k\ y =0 } \textbf{ for all } \pmb{k\ge 0} \textbf{ and } \pmb{m_{n}=0}\hspace{3cm}		\hfill {\pmb{(\mathcal{SM})}}
\]
\end{criterion}
Note that the matrices $(P_{n-2}X^k)$, $k\ge 1$, are symmetric.

\begin{proof}
We start from $(\mathcal{ S})$ and remember that $wQ^{-1}v=c_{n+1}$ just means $m_n=\det M_0=0$. 
For $k\ge 1$ one has $wQ^{-1}(BQ^{-1})^kv=wQ^{-1}(BQ^{-1})^{k-1}BQ^{-1}v$, where $(BQ^{-1})^0=I_{n-1}$, $BQ^{-1}=\bmx X&\star\\0&0\emx$, $BQ^{-1}v=B(Q^{-1}v)=\bmx y\\0\emx$, $wQ^{-1}=\bmx \tr yP_{n-2}&m_1\emx$.
For $k=1$ this gives 
$wQ^{-1}(BQ^{-1}v)=\bmx \tr yP_{n-2}&m_1\,\emx\bmx y\\0 \emx=\tr yP_{n-2}\,y$. 
For \mbox{$k>1$} we have $(BQ^{-1})^{k-1}=\bmx X^{k-1}&\star\\0&0\emx$ and therefore $(BQ^{-1})^{k-1}(BQ^{-1})v=\bmx X^{k-1}y\\0\emx$ and
$wQ^{-1} (BQ^{-1})^{k-1}(BQ^{-1})v=\tr yP_{n-2}X^{k-1}y$. This gives us $(\mathcal{SM})$. Similarly one can translate back to $(\mathcal S)$.
\end{proof}

As a result ToePC becomes equivalent to the conjecture
\begin{center}
{\bf 
$\pmb{(\mathcal{SM})}$ implies $\pmb{y=0}$
}
\end{center}
This new algebraic conjecture looks quite different compared with ToePC and proposes different methods to be proved. Simulations with computer algebra software show that considering the minors as independent variables (compare (\ref{var})) leads to equations that only admit the zero solution and thusly confirm the conjecture.

To conclude I have to thank Alexander Kovacec for several not only mathematical discussions and encouragement.


\begin{thebibliography}{al}
\bibitem{ScSh1}
W. Schmale, P. K. Sharma: \textit{Cyclizable matrix pairs over C[x] and a conjecture on T\"oplitz matrices}, Linear Algebra and its Applications, {389} (2004), pp. 33-42.
\bibitem{BSSV} R. Bumby, E. Sontag, H. Sussmann, W. Vasconcelos: \textit{Remarks
on the pole shifting problem over rings}, J. Pure Appl. Algebra 20
(1981), pp. 113-127.
\bibitem{ScSh2}
W. Schmale, P. K. Sharma: \textit{Problem 30-3: Singularity of a Toeplitz Matrix}, IMAGE {issues Nr. 30 to 37} (2003-2006).
\bibitem{BoLi}V. Bolotnikov, Ch.-K. Li: \textit{Solution to Problem 30-3}, IMAGE {36} (2006), p.26.
\bibitem{HW} H. Wimmer: \textit{Remark on Problem 30-3}, IMAGE {37} (2006), p. 19.
\bibitem{open}
W. Schmale: \textit{The T\"oplitz Pencil Conjecture is still open, Letter to the editors of Linear and multilinear algebra}, Linear and multilinear algebra, {63, No 3} (2015), p. 650.
\bibitem{GoKo}
A. Kovacec, M. C. Gouveia: \textit{The Hankel pencil conjecture},
Linear Algebra and its Applications
{431} (2009), pp. 1509-1525.
\bibitem{G} F.R. Gantmacher: \textit{The theory of matrices}, Chelsea Publ. (1974), Vol. II.
\bibitem{T} W. F. Trench: \textit{Inverses of lower triangular Toeplitz matrices}, (2009), unpublished note available at: 
\verb| http://works.bepress.com/william_trench/132/ |
\bibitem{HoJo} R.A. Horn, C.R. Johnson: \textit{Matrix Analysis, Second Edition}, Cambridge University Press (2013).


\end{thebibliography}
\end{document}